\documentclass[11pt]{article}
\usepackage{}
\usepackage[total={6in, 9in}]{geometry}
\usepackage{amsfonts}
\usepackage{amsthm}
\usepackage{amssymb}
\usepackage{amsmath,amsfonts}
\usepackage{mathrsfs}
\usepackage{cases}
\usepackage{latexsym,bm}
\usepackage{indentfirst}
\usepackage{color}
\usepackage{ifpdf}
\usepackage{graphicx}
\usepackage{psfrag}
\usepackage{enumerate}
\usepackage{dsfont}
\usepackage[
pdfauthor={SS},
pdftitle={Hamiltonian cyccles in $4$-tough $(P_2\cup kP_1)$-free graphs},
pdfstartview=XYZ,
bookmarks=true,
colorlinks=true,
linkcolor=blue,
urlcolor=blue,
citecolor=blue,
bookmarks=true,
linktocpage=true,
hyperindex=true
]{hyperref}
\newtheorem{thm}{Theorem}
\newtheorem{cor}[thm]{Corollary}
\newtheorem{lem}[thm]{Lemma}

\newtheorem{conj}[thm]{Conjecture}

\newtheorem{CLA}{\textbf{Claim}}
\newcommand{\iC}{\overset{\leftharpoonup }{C}}
\newcommand{\oC}{\overset{\rightharpoonup }{C}}
\begin{document}
	
	\title{A note on hamiltonian cycles in $4$-tough $(P_2\cup kP_1)$-free graphs}
\author{Lingjuan Shi$^{a}$\footnote{
The first author is supported by NSFC (grant no. 11871256 and 11901458) and by the Fundamental Research Funds for the Central Universities (grant no. D5000200199).}, Songling Shan$^b$
\date{\small $^{a}$School of Software, Northwestern Polytechnical University, Xi'an, Shaanxi 710072, P. R. China\\
$^{b}$ Department of Mathematics, Illinois State University, Normal, IL 61790, USA\\
E-mails: shilj18@nwpu.edu.cn, sshan12@ilstu.edu}
}
\maketitle
\begin{abstract}
Let $t>0$ be a real number and $G$ be a graph. We say $G$
is $t$-tough if for every cutset $S$ of $G$, the ratio of
$|S|$ to the number of components of $G-S$ is at least $t$.  The Toughness Conjecture of Chv\'atal, stating that there exists a constant $t_0$ such that every $t_0$-tough graph with at least three vertices  is hamiltonian, is still open in general.
For any given  integer $k\ge 1$, a graph $G$ is $(P_2\cup kP_1)$ free if $G$ does not
contain the disjoint union of $P_2$ and $k$ isolated vertices as an induced subgraph.
 In this note, we show that every 4-tough and $2k$-connected  $(P_2\cup kP_1)$-free graph with at least three vertices is hamiltonian.
 This result in some sense is an ``extension'' of the classical Chv\'{a}tal-Erd\H{o}s Theorem that every $\max\{2,k\}$-connected
 $(k+1)P_1$-free graph on at least three vertices is hamiltonian.

\medskip


\end{abstract}

\section{Introduction}
We  consider only  simple graphs. Let $G$ be a graph. We denote  by $V(G)$ and $E(G)$ the vertex set and edge set of $G$, respectively.
For a vertex $x\in V(G)$ and a subgraph $H$ of $G$,  $N_H(x)$
denotes the set  of neighbors of $x$ that are contained in $V(H)$,  and $d_H(x):=|N_H(x)|$.
Let  $S\subseteq V(G)$. Then $N_G(S)=\bigcup_{x\in S}N_G(x)\setminus S$ and $N_H(S)=N_G(S)\cap V(H)$.
 Denote by $G[S]$  the  subgraph of $G$ induced by $S$,  and  let $G-S=G[V(G)\setminus S]$.
We skip the  subscript $G$  if no confusion may arise.

For a given graph $H$, we say that $G$ is $H$-free if $G$ does not
contain  $H$ as an induced subgraph.  For another graph $F$, $H\cup F$
is the disjoint union of $H$ and $F$.
For positive  integers $a$ and $b$, we denote by $aP_b$ the graph consisting of $a$ disjoint copies of the path $P_b$.
For two integers $p$ and $q$, let $[p,q]=\{i\in \mathbb{Z}: p\le  i \le q\}$.

Throughout this paper,  if not specified,
we will assume $t$ to be a nonnegative real number.
The number of components of $G$ is denoted by $c(G)$.  The graph $G$ is said to be  \emph{$t$-tough} if $|S|\ge t\cdot
c(G-S)$ for each $S\subseteq V(G)$ with $c(G-S)\ge 2$. The \emph {toughness $\tau(G)$} is the largest real number $t$ for which $G$ is
$t$-tough, or is  $\infty$ if $G$ is complete. This concept, a
measure of graph connectivity and ``resilience'' under removal of
vertices, was introduced by Chv\'atal~\cite{chvatal-tough-c} in 1973.
It is  easy to see that  if $G$ has a hamiltonian cycle
then $G$ is 1-tough. Conversely,
Chv\'atal~\cite{chvatal-tough-c}
conjectured that
there exists a constant $t_0$ such that every
$t_0$-tough graph is hamiltonian (Chv\'atal's toughness conjecture).
Bauer, Broersma and Veldman~\cite{Tough-counterE} have constructed
$t$-tough graphs that are not hamiltonian for all $t < \frac{9}{4}$, so
$t_0$ must be at least $\frac{9}{4}$. It is not difficult to see that
a non-complete $t$-tough graph is $2\lceil t\rceil$-connected.

There are many  papers on
Chv\'atal's toughness conjecture,
and it  has
been verified when restricted to a number of graph
classes~\cite{Bauer2006},
including planar graphs, claw-free graphs, co-comparability graphs, and
chordal graphs. The conjecture was also confirmed for
graphs forbidden small linear forests.   The results include the following:  1-tough
for $R$-free
graphs, where  $R\in\{P_4,  P_3\cup P_1, P_2\cup 2P_1\}$~\cite{MR3557210};
25-tough, 3-tough, and now 2-tough for $2P_2$-free graphs~\cite{2k2-tough, MR4102606, 2103.06760};
3-tough for $(P_2\cup3P_1)$-free graphs~\cite{2106.07083};
7-tough for $(P_3\cup2P_1)$-free graphs~\cite{2107.08476}; and
 15-tough  for $(P_2\cup P_3 )$-free graphs~\cite{MR4245269}.
 In this paper, for any integer $k\ge 1$, we  investigate Chv\'atal's toughness conjecture for $(P_2\cup kP_1)$-free graphs.  For $k=1$, it is an easy exercise to show that every 1-tough $(P_2\cup P_1)$-free graph on $n\ge 3$ vertices has minimum degree at least $n/2$ and so is hamiltonian by Dirac's Theorem~\cite{MR47308}.  By~\cite{MR3557210}, 1-tough is enough to guarantee a hamiltonian cycle
for $(P_2\cup 2P_1)$-free graphs on at least three vertices. By~\cite{2106.07083}, 3-tough is enough to guarantee a hamiltonian cycle
for $(P_2\cup 3P_1)$-free graphs on at least three vertices. Thus we consider only the case when $k \ge 4$, and obtain the result below.


\begin{thm}\label{main}
Let $k\geq 4$ be an integer and $G$ be a $4$-tough and $2k$-connected $(P_2\cup kP_1)$-free graph. Then $G$ is hamiltonian.
\end{thm}

As a non-complete $t$-tough graph is $2\lceil t\rceil$-connected, as a corollary of Theorem~\ref{main},
we get the following result.

\begin{cor}\label{cor}
	For any integer $k\ge 4$, every $k$-tough $(P_2\cup kP_1)$-free graph on at least three vertices is hamiltonian.
\end{cor}

Theorem~\ref{main} can also be seen as a result along the Chav\'{a}tal-Erd\H{o}s type theorems~\cite{MR297600}.
The Chav\'{a}tal-Erd\H{o}s Theorem can be rephrased as below:  for any integer $k\ge 1$, every $\max\{2,k\}$-connected
$(k+1)P_1$-free graph on at least three vertices is hamiltonian.  It is necessary  for a hamiltonian graph to be 1-tough and the condition that  ``$k$-connected and $(k+1)P_1$-free'' in the  Chav\'{a}tal-Erd\H{o}s Theorem implies that the graph is 1-tough.
However, any constant connectivity condition cannot guarantee  the existence of a hamiltonian cycle in
 $(P_2\cup kP_1)$-free graphs.  For example,  for  integer $n\ge 2$, the complete bipartite graph $K_{n-1,n}$  is  $(n-1)$-connected, $(P_2\cup kP_1)$-free for any $k\ge 1$, but is not hamiltonian.
 This explains the involvement of the toughness condition in Theorem~\ref{main} other than the connectivity condition. However, we do think that the condition 4-tough is not sharp and we propose the following conjecture.

\begin{conj}
	Let $k\geq 4$ be an integer and $G$ be a $1$-tough and $2k$-connected $(P_2\cup kP_1)$-free graph. Then $G$ is hamiltonian.
	\end{conj}
\section{Proof of Theorem \ref{main}}
We start with certain notation and preliminary results.
Let $C$ be an oriented cycle, and we assume that the orientation is clockwise throughout the rest of this paper. For $x\in V(C)$,
denote the immediate successor of $x$ on $C$ by $x^+$ and the immediate  predecessor of $x$ on $C$ by $x^-$.
For $u,v\in V(C)$, $u\oC v$  denotes the segment of $C$
starting at $u$, following $C$ in the orientation,  and ending at $v$.
Likewise, $u\iC v$ is the opposite segment of $C$ with endpoints as $u$
and $v$.

A path $P$ connecting two vertices $u$ and $v$ is called
a {\it $(u,v)$-path}, and we write $uPv$ or $vPu$ in order to specify the two endvertices of
$P$. Let $uPv$ and $xQy$ be two paths. If $vx$ is an edge,
we write $uPvxQy$ as
the concatenation of $P$ and $Q$ through the edge $vx$.

\begin{lem}[\cite{MR1336668}]\label{Bauer}
Let $t$ be a positive real number and $G$ be a $t$-tough graph on at lest three vertices
with $\delta(G)>\frac{n}{t+1}-1$. Then $G$ is hamiltonian.
\end{lem}

\proof[\bf Proof of Theorem~\ref{main}]
Suppose for a contradiction that $G$ is not hamiltonian. Then $G$ is not complete and  $\delta(G)\leq\frac{n}{5}-1$ by  Lemma~\ref{Bauer}. Since $G$ is   $2k$-connected,  $\delta(G)\geq 2k$.
Thus we have  $n\geq5(2k+1)\geq 45$.
Let $C$ be a longest cycle of $G$. As $G$ is not hamiltonian, $V(G)\setminus V(C)\neq\emptyset$.
Let $H$ be a component of $G-V(C)$ and
suppose $|N_C(V(H))|=t$. We let $x_1, \ldots, x_t$ be all the neighbors
of vertices of $H$ on $C$, and we assume that they appear in the order $x_1, \ldots, x_t$ along $\oC$.
Note that $t \ge 2k$ as
 $G$ is $2k$-connected.

\begin{CLA}\label{claim:neighbors}
Any two vertices in $\{x_1,\ldots, x_t\}$ are not consecutive on $C$,
and $\{x_1^+, x_2^+, \ldots, x_t^+\}$ is an independent set in $G$.
\end{CLA}

\proof The first part is clear as otherwise we can extend $C$.
For the second part,
suppose  to the contrary that  $x_i^+x_j^+\in E(G)$ for some $i,j\in[1,t]$.
Without loss of generality, assume that $i<j$.
Let $x_i', x_j'\in V(H)$ such that $x_ix_i', x_jx_j'\in E(G)$,
and let $P$ be an $(x_i', x_j')$-path of $H$.
Then $x_ix_i'Px_j'x_j \iC x_i^+x_j^+ \oC x_i$  is a cycle longer than $C$, a contradiction.
\qed

\begin{CLA}\label{claim:trivial-components}
	Each component of $G-V(C)$ is trivial, i.e., a one-vertex graph.
\end{CLA}

\proof
We still let $H$ represent any component of $G-V(C)$.
Suppose to the contrary that
$H$ has an edge $uv$.  Then $uv$ and  $\{x_1^+, x_2^+, \ldots, x_t^+\}$  together form an induced  $P_2\cup tP_1$, which contains $P_2\cup kP_1$ as an induced subgraph, a contradiction.
\qed

\begin{CLA}\label{claim:C-length}
	$|V(C)|\geq\frac{4n}{5}$.
\end{CLA}

\proof
Suppose to the contrary that $|V(C)|<\frac{4n}{5}$. Then $G-V(C)$
has exactly  $n-|V(C)|>n/5\ge 9$ components by Claim~\ref{claim:trivial-components}.
Thus $V(C)$ is a cutset of $G$ and $\frac{|V(C)|}{c(G-V(C))}<4$.
This contradicts that $G$ is $4$-tough.
\qed

Our goal in the rest of the proof  is to find an independent set of size more than $n/5$
in $G$ and so to reach a contradiction to the toughness of $G$.
By Claim~\ref{claim:trivial-components}, let $V(H)=\{x\}$.
Assume, without loss of generality,  that  $|V(x_1 \oC x_{k+1})|\leq | V(x_{k+1}\oC x_1)|$.
Then $| V(x_{k+1}\oC x_1)| \ge \frac{2n}{5}$ by Claim~\ref{claim:C-length}.
Let
 $x_{k+1} \oC x_1=x_{k+1}x_{k+1}^+y_1y_2\ldots y_h x_1$, and
 define
 \begin{numcases} {Y=}
 	\{y_2, y_4, \ldots, y_h\}  & \text{if $h$ is even};  \nonumber \\
 	\{y_2, y_4, \ldots, y_{h-1}\} & \text{if $h$ is odd}.  \nonumber
 \end{numcases}

\begin{CLA}\label{claim:ind-Y}
$N(Y)\cap \{x_1^+, \ldots, x_{k+1}^+\}=\emptyset$.
As a consequence, $Y$ is an independent set in $G$.
\end{CLA}
\proof

Let  $X=\{x_1^+, x_2^+, \ldots, x_{k+1}^+\}$. Since $X$ is an independent set
in $G$ by Claim~\ref{claim:neighbors}, the consequence part of the claim
follows easily from  the first part and the $(P_2\cup kP_1)$-freeness of $G$.
So we only show that $N(Y)\cap \{x_1^+, \ldots, x_{k+1}^+\}=\emptyset$.
Considering  the edge  $x_{k+1}^+y_1$ and the independent set $X$,
we conclude that $|N(y_1)\cap X|\geq2$ by the
$(P_2\cup kP_1)$-freeness of $G$.
Next, we show $N(y_2)\cap X=\emptyset$. Suppose not, then $| N(y_2)\cap X|\geq2$.
Otherwise let $ N(y_2)\cap X=\{x_i\}$ for some $i\in [1,k+1]$. Then $y_2x_i$
and $X\setminus\{x_i\}$ form  an induced copy of $(P_2\cup kP_1)$.
We consider three cases regarding the size of the set  $N(y_1)\cap N(y_2)\cap X$
below.

{  \noindent \textbf Case 1:} $|N(y_1)\cap N(y_2)\cap X|\geq2$.

Let  $x_i^+, x_j^+\in N(y_1)\cap N(y_2)\cap X$  for some distinct $i,j\in [1,k+1]$. Assume, without loss of generality, that $i<j$. Then $xx_i\iC y_2x_j^+ \oC y_1x_i^+ \oC x_jx$ is a cycle longer than $C$, a contradiction.

{ \noindent \textbf Case 2:} $|N(y_1)\cap N(y_2)\cap X|=1$.

Let $N(y_1)\cap N(y_2)\cap X=\{x_i^+\}$ for some $i\in [1,k+1]$.
If $i=1$, then $y_2$ has another neighbor $x_j^+$ with $ j\in [2,k+1]$ since $|N(y_2)\cap X|\geq2$. So $xx_1\iC y_2x_j^+\oC y_1x_1^+\oC x_jx$ is a cycle longer than $C$, a contradiction.
If $i=k+1$, then $y_1$ has another neighbor $x_j^+$ with $j\in [1,k]$ since $|N(y_1)\cap X|\geq2$.
So $xx_j\iC y_2x_{k+1}^+\oC y_1x_j^+\oC x_{k+1}x$ is a cycle longer than $C$, a contradiction.
If $i\in[2,k]$, then we define two indices as follows:
\begin{eqnarray*}
s&=&\min\{j: x_j^+\in N(y_1)\cap X\},\\
\ell&=&\max\{j: x_j^+\in N(y_2)\cap X\}.
\end{eqnarray*}
Clearly, $s\leq i$ and $\ell \geq i$. If $s<i$ then $xx_{s}\iC y_2x_i^+\oC y_1x_{s}^+\oC x_ix$ is a cycle longer than $C$; if
 $\ell >i$, then $xx_i\iC y_2x_{\ell }^+\oC y_1x_i^+\oC x_{\ell }x$  is a cycle longer than $C$.
Thus we assume $s=\ell=i$.

If $y_1$ is a neighbor of $x$, then $xy_1\iC x_i^+y_2\oC x_ix$ is a cycle  longer than $C$, a contradiction.
So $y_1\notin N(x)$. This implies that $y_2\notin\{x_1^+, x_2^+, \ldots, x_t^+\}$.
We claim that $N(y_1)\cap\{x_{k+2}^+, x_{k+3}^+, \ldots, x_t^+\}=\emptyset$.  Otherwise,
suppose  $y_1x_h^+\in E(G)$ for some $h\in [k+2, t]$.
Then $xx_i\iC x_h^+y_1 \iC x_i^+y_2 \oC x_hx$
is  a cycle longer than $C$.
Hence  $x_i^+\oC y_1$ or $y_2\oC x_i^+$ has at least $k+1$ vertices belonging to $\{x_1^+, x_2^+, \ldots, x_t^+\}$ as $t\ge 2k$.

If $x_i^+\oC y_1$ has at least $k+1$ vertices in $\{x_1^+, x_2^+, \ldots, x_t^+\}$, then $y_2x_i^+$ and the $k$ vertices in $(V(x_i^+\oC y_1)\cap\{x_1^+, x_2^+, \ldots, x_t^+\})\setminus \{x_i^+\}$ induce a $P_2\cup kP_1$ subgraph in $G$, a contradiction.
If $y_2\oC x_i^+$ has at least $k+1$ vertices in $\{x_1^+, x_2^+, \ldots, x_t^+\}$, then $x_i^+y_1$ and the $k$ vertices in $(V(y_2\oC x_i^+)\cap\{x_1^+, x_2^+, \ldots, x_t^+\})\setminus \{x_i^+\}$ induce a $P_2\cup kP_1$ subgraph in $G$, a contradiction.

{ \noindent \textbf Case 3:} $|N(y_1)\cap N(y_2)\cap X|=0$.

We define $s$ and $\ell $ the same way as in Case 2. Since $|N(y_1)\cap N(y_2)\cap X|=0$, $s\ne \ell$.
If $s<\ell $, then $xx_{s}\iC y_2x_{\ell }^+\oC y_1x_{s}^+\oC x_{\ell }x$ is a cycle longer than $C$, a contradiction.
So $s>\ell $.
If $y_1\in N(x)$,  then  $xx_{\ell }\iC y_2x_{\ell }^+\oC y_1x$  is a cycle longer than $C$, a contradiction.
So $y_1$ is not a neighbor of $x$ and  thus $y_2\notin\{x_1^+, x_2^+, \ldots, x_t^+\}$.
By the same argument as in Case 2, we have $N(y_1)\cap\{x_{k+2}^+, x_{k+3}^+, \ldots, x_t^+\}=\emptyset$.
So $y_1$ dose not have  any neighbor  in $\{x_1^+, \ldots, x_{s-1}^+\} \cup \{x_{k+2}^+,  \ldots, x_t^+\}$, which is a subset of $V(y_2\oC x_{s})$.  By the definition of $\ell$ and the assumption that $\ell <s$, we know that $y_2$ does not have any neighbor in $\{x_s^+, \ldots, x_{k+1}^+\}$.
Since $t\geq2k$, $y_2\oC x_{s}$ or $x_{s}^+\oC x_{k+1}^+$ has at least $k$ vertices belonging to $\{x_1^+, x_2^+, \ldots, x_t^+\}$.
If $y_2\oC x_{s}$ has $k$ vertices belonging to $\{x_1^+, x_2^+, \ldots, x_t^+\}$, then those $k$ vertices and the edge $y_1x_{s}^+$ induce a copy of $P_2\cup kP_1$.
If $x_{s}^+\oC x_{k+1}^+$ has $k$ vertices belonging to $\{x_1^+, x_2^+, \ldots, x_t^+\}$, then those $k$ vertices and  the edge $y_2x_{\ell }^+$ induce a copy of $P_2\cup kP_1$.

 Cases~1 to 3  imply $N(y_2)\cap X=\emptyset$.
Since $X$ is an independent set in $G$, $|X|=k+1$ and $y_2y_3\in E(G)$,  the $(P_2\cup kP_1)$-freeness of $G$ implies that $|N(y_3)\cap X|\geq2$.  With $y_3$ playing the role of
$y_1$ and $y_4$ playing the role of $y_2$, the same arguments in Cases~1 to 3
show that $N(y_4)\cap X=\emptyset$.  Repeating the  same arguments for all  other elements of
$Y$, we get $N(Y)\cap X=\emptyset$.
\qed

By Claim~\ref{claim:neighbors} and Claim~\ref{claim:ind-Y}, we know that  $X\cup Y$ in an independent set in $G$. By Claim~\ref{claim:C-length}  and the assumption  that $|V(x_1\oC x_{k+1})|\leq|V(x_{k+1}\oC x_1)|$, we have $|Y|\geq\lfloor\frac{1}{2}(\frac{|C|-2}{2})\rfloor\geq\lfloor\frac{n}{5}-\frac{1}{2}\rfloor$.
Since $k\ge 4$ and so $|X| \ge 5$,
we have  $|X\cup Y|>\frac{n}{5} \ge 9$.  Let $S=V(G)\setminus (X\cup Y)$.
Then $S$ is a cutset of $G$ with $c(G-S)=|X\cup Y|>\frac{n}{5}$.
However, $\frac{|S|}{c(G-S)}<4$, contradicting $G$ being 4-tough.
The proof is now complete.
\qed


\bibliographystyle{plain}
\bibliography{SSL-BIB}

\end{document}